\documentclass[a4paper, twoside, psamsfonts]{amsart}

\usepackage{graphicx}
\usepackage{subfigure}
\usepackage{booktabs}

\numberwithin{equation}{section}

\newcommand{\velocity}{u}
\newcommand{\testvector}{v}
\newcommand{\pressure}{p}
\newcommand{\testscalar}{q}

\newcommand{\R}{\mathbb{R}}
\newcommand{\V}{\mathbb{V}}

\newcommand{\VV}{V}
\newcommand{\QQ}{Q}

\newcommand{\Ltwo}{L^2(\Omega)}
\newcommand{\Hdiv}{H(\Div, \Omega)}
\newcommand{\LtwoVec}{L^2(\Omega; \V)}

\newcommand{\cP}[1]{\mathcal{P}^{c}_{#1}}
\newcommand{\dP}[1]{\mathcal{P}_{#1}}
\newcommand{\Poly}[1]{\mathcal{P}_{#1}}

\DeclareMathOperator{\Div}{div}

\DeclareMathOperator{\Grad}{grad}

\newcommand{\bigO}{\mathcal{O}}
\newcommand{\triang}{\mathcal{T}}
\newcommand{\foralls}{\forall \,}

\newcommand{\contBabuska}{\gamma}
\newcommand{\contCoercivity}{\alpha}
\newcommand{\contInfsup}{\beta}
\newcommand{\Babuska}{\gamma_h}
\newcommand{\Coercivity}{\alpha_h}
\newcommand{\Infsup}{\beta_h}
\newcommand{\contLaplaceInfsup}{\beta^{\Div}}
\newcommand{\LaplaceInfsup}{\beta_h^{\Div}}
\newcommand{\StokesInfsup}{\beta_h^{1}}
\newcommand{\ReducedInfsup}{\tilde{\beta}_h}

\newcommand{\digit}{\hphantom0}

\newcommand{\Kernel}{\mathcal{Z}}
\newcommand{\Spurious}{\mathcal{N}}

\theoremstyle{plain}
\newtheorem{thm}{Theorem}[section] 
\newtheorem{lem}[thm]{Lemma}

\newtheorem{definition}[thm]{Definition}

\theoremstyle{remark}
\newtheorem*{remark}{Remark}
\begin{document} 

\title{Stability of Lagrange elements for the mixed Laplacian}

\author{Douglas N. Arnold}
\address{School of Mathematics,
512 Vincent Hall, 206 Church Street S.E.,
University of Minnesota,
Minneapolis, MN 55455}
\email{arnold@umn.edu}
\thanks{The work of the first author was supported in part by NSF grant
DMS-0713568.}

\author{Marie E. Rognes}
\address{Center for Biomedical Computing, 
  Simula Research Laboratory, P.O.Box 134, 1325 Lysaker, Norway.
  }
\email{meg@simula.no}
\thanks{The work of the second author was supported by a Center of Excellence grant from the
  Norwegian Research Council to Centre of Mathematics for
  Applications at the University of Oslo and by a Center of Excellence grant from the Norwegian
  Research Council to Center for Biomedical Computing at Simula
  Research Laboratory.}

\subjclass[2000]{Primary: 65N30}
\keywords{mixed finite elements, Lagrange finite elements, stability}

\begin{abstract}
The stability properties of simple element choices for the mixed
formulation of the Laplacian are investigated numerically.  The element
choices studied use vector Lagrange elements, i.e., the space of
continuous piecewise polynomial vector fields of degree at most $r$,
for the vector variable, and the divergence of this space, which consists of
discontinuous piecewise polynomials of one degree lower, for the scalar
variable.  For polynomial degrees $r$ equal 2 or 3, this pair of spaces
was found to be stable for all mesh families tested. In particular,
it is stable on diagonal mesh families, in contrast to its behaviour
for the Stokes equations. For degree $r$ equal 1, stability holds for
some meshes, but not for others. Additionally, convergence was observed
precisely for the methods that were observed to be stable. However,
it seems that optimal order $L^2$ estimates for the vector variable,
known to hold for $r>3$, do not hold for lower degrees.
\end{abstract}

\maketitle

\section{Introduction}


In this note, we consider approximations of the mixed Laplace
equations with Dirichlet boundary conditions: Given a source $g$, find
the velocity $\velocity$ and the pressure $\pressure$ such that
\begin{equation*}
  \velocity - \Grad p = 0, \; \Div \velocity = g \; \text{ in } \Omega,
  \quad \quad
  p = 0  \; \text{ on } \partial \Omega,
\end{equation*}
for a domain $\Omega \subset \R^2$ with boundary $\partial \Omega$.
The equations offer the classical weak formulation: Find a square
integrable vector field with square integrable divergence $\velocity
\in H(\Div, \Omega)$ and a square integrable function $p \in
L^2(\Omega)$ such that
\begin{equation}
  \label{eq:laplace}
  \int_{\Omega} \velocity \cdot \testvector
  +  \testscalar \, \Div \velocity 
  + \pressure \,  \Div \testvector
  = \int_{\Omega} g \, \testscalar
\end{equation}
for all $\testvector \in \Hdiv$ and $\testscalar \in \Ltwo$.  The
above formulation can be discretized using a pair of finite
dimensional spaces $V_h \subset H(\Div,\Omega)$, $Q_h\subset L^2(\Omega)$, yielding discrete approximations
$\velocity_h \in \VV_h$ and $\pressure_h \in \QQ_h$
satisfying~\eqref{eq:laplace} for all $\testvector \in \VV_h$ and
$\testscalar \in \QQ_h$.


As is well-known, the spaces $V_h$ and $Q_h$ must satisfy certain
stability, or compatibility, conditions for the discretization to be
well-behaved~\cite{art:Bre1974}. More precisely, there must exist
positive constants $\contCoercivity$ and $\contInfsup$ such that for
any $h$,
\begin{subequations}
  \begin{align}
    \label{eq:laplace:coercivity}
    0 < \contCoercivity \leq \Coercivity = 
    \inf_{\velocity \in \Kernel_h} \sup_{\testvector \in \Kernel_h}
    \frac{\langle \velocity, \testvector \rangle}
         {\|\velocity\|_{\Div} \|\testvector\|_{\Div}}, \\
    \label{eq:laplace:infsup}
    0 < \contInfsup < \LaplaceInfsup 
    = 
    \inf_{\testscalar \in \QQ_h} \sup_{\testvector \in \VV_h}
    \frac{\langle \Div \testvector, \testscalar  \rangle}
         {\|\testvector\|_{\Div} \|\testscalar\|_0}.
  \end{align}
\end{subequations}
Here, $\|\cdot\|_{\Div}$ and $\|\cdot\|_{0}$ denote the norms on
$\Hdiv$ and $\Ltwo$, $\langle \, \cdot \, , \, \cdot \, \rangle$ is
the $\Ltwo$ inner product and
\begin{equation}
  \Kernel_h = \{\testvector \in \VV_h \, | \; 
  \langle \Div \testvector, \testscalar \rangle
  = 0 \quad \foralls \testscalar \in \QQ_h \}.
\end{equation}
The two conditions will be referred to as the Brezzi coercivity and
the Brezzi inf-sup condition for the mixed Laplacian. The classical
conforming discretizations of~\eqref{eq:laplace} rely on the finite
element families of Raviart and Thomas~\cite{art:RavTho1977} or
Brezzi, Douglas and Marini~\cite{art:BreDouMar1985} for the space $V_h
\subset H(\Div)$ in order to satisfy these conditions.


In this note, we shall consider the Lagrange vector element spaces,
that is, continuous piecewise polynomial vector fields defined
relative to a triangulation $\triang_h$, for the space $V_h$. This is
motivated by the following reasons. First, these spaces are fairly
inexpensive, simple to implement and post-process and in frequent use
for other purposes. Second, such pairs would allow continuous
approximations of the velocity variable, or when viewed in connection
with linear elasticity, lay the ground for continuous approximations
of the stress tensor. Moreover, in the recent years, there has been an
interest in mixed finite element discretizations that are both stable
for~\eqref{eq:laplace} and for the Stokes equations:
\begin{equation}
  \label{eq:stokes}
  \int_{\Omega} \Grad \velocity : \Grad \testvector
  + \testscalar \, \Div \velocity
  + \pressure \, \Div \testvector
  = \int_{\Omega} f \, \testvector
\end{equation}
for all $\testvector \in H^1(\Omega; \V)$ such that $\int_{\Omega} v =
0$ and all $\testscalar \in \Ltwo$.  The search for conforming such
discretizations is complicated by the fact that the existing, stable
discretizations of~\eqref{eq:laplace} are such that $V_h \not \subset
H^1(\Omega; \V)$. On the other hand, the existing stable
discretizations of~\eqref{eq:stokes} are typically unstable
for~\eqref{eq:laplace}~\cite{art:MarTaiWin2002}. The existence of
stable discretizations $V_h \times Q_h$ of~\eqref{eq:laplace}
such that $V_h \subset H^1(\Omega; \V)$ becomes a natural separate
question. Unfortunately, there are no known such finite element
discretizations that are stable for any admissible triangulation
family $\{\triang_h\}$. In this note, we aim to numerically examine
cases where a reduced stability property may be identified. In this
sense, the investigations here are in the spirit of the work of
Chapelle and Bathe~\cite{art:ChaBat1993} and
Qin~\cite{thesis:Qin1994}.


For a family of conforming discretizations $\{\VV_h \times
\QQ_h\}_{h}$ of~\eqref{eq:laplace} such that $\Div \VV_h \subseteq
\QQ_h$ for each $h$, the condition~\eqref{eq:laplace:coercivity} is
trivial. The stability conditions thus reduce to the
condition~\eqref{eq:laplace:infsup}, namely the question of bounded
Brezzi inf-sup constant $\LaplaceInfsup$. On the other hand, recall
that for the Stokes formulation~\eqref{eq:stokes}, the corresponding
Brezzi coercivity condition is trivial by the Poincar\'e inequality. Hence,
for any family of conforming discretizations, the stability conditions
for Stokes reduce to that of a uniform bound for the Brezzi inf-sup
constant $\StokesInfsup$. Here,
\begin{equation}
  \label{eq:H1_Hdiv_inequality}
  \StokesInfsup = \inf_{\testscalar \in \QQ_h} \sup_{\testvector \in
  \VV_h} \frac{\langle \Div \testvector, \testscalar
  \rangle}{\|\testvector\|_{1} \|\testscalar\|_0}.
\end{equation}
when $\VV_h \subset H^1(\Omega; \V)$, $\QQ_h \subset L^2(\Omega)$ and
$\|\cdot\|_{1}$ denotes the norm on $H^1(\Omega)$.  Further, such a
bound immediately gives~\eqref{eq:laplace:infsup} since $\StokesInfsup
\leq \LaplaceInfsup$ by definition. Hence, if $\Div \VV_h \subseteq
\QQ_h$, stability for Stokes immediately gives stability for the mixed
Laplacian.

The conditions $\Div \VV_h \subseteq \QQ_h$ and $\VV_h \subset
H^1(\Omega; \V)$ are clearly satisfied by the element pairs consisting
of continuous piecewise polynomial vector fields of degree less than
or equal to $r$ and discontinuous piecewise polynomials of degree
$r-1$, for $r = 1, 2, \dots$.  This family could be viewed as an
attractive family of elements for both the Stokes equations and the
mixed Laplacian. However, the Brezzi inf-sup constant(s) will not be
bounded for all $r$. For $r \geq 4$, Scott and Vogelius demonstrated
that these finite element spaces will be stable for the Stokes
equations on triangulations that have no nearly singular vertices,
that is, triangulations that are not singular in the appropriate
sense~\cite{art:ScoVog1985}. The lower order cases, $1 \leq r \leq 3$,
were studied carefully by Qin, concluding that the elements are not
stable in general~\cite{thesis:Qin1994}. However, they are stable for
some specific families of triangulations, and can be stabilized by
removal of spurious pressure modes on some other classes of
triangulations. (The space $\Spurious_h$ of spurious pressure modes is
defined in~\eqref{eq:kernelcomplement} below.) In general, the
stability of finite element spaces for the Stokes equations has been
extensively investigated. In addition to the previous references,
surveys are presented in~\cite{art:BofBreGas2007,
book:BreFor1991}. However, to our knowledge, a careful study of the
lower order cases has not been conducted for the mixed Laplacian.

As the stability for the mixed Laplacian is a weaker requirement when
$\Div \VV_h \subseteq \QQ_h$, there may be a greater class of
triangulations for which the elements form a stable discretization. In
fact, this is known to be true. One example is provided by the pairing
of continuous piecewise linear vector fields $V_h$ and the subspace of
discontinuous piecewise constants $Q_h$ such that $Q_h = \Div V_h$ on
crisscross triangulations of the unit square. Qin proved that there
does not exist a $\contInfsup > 0$ such that $\StokesInfsup >
\contInfsup$ for any $h$~\cite[Lemma 7.3.2]{thesis:Qin1994}. On the
other hand, Boffi et al.~proved that such a bound does exists for
$\LaplaceInfsup$~\cite{art:BofBreGas1999}. We shall present numerical
evidence suggesting that there is a range of triangulations for which
this phenomenon occurs. The main results are summarized below.
\begin{description}
  \item[Spurious modes] For $r = 2, 3$ and for all triangulations
  tested, the dimension of the space $\Spurious_h$ of spurious modes
  is equal to the number of interior singular vertices
  $\sigma$. However, for $r = 1$ and one of the triangulation families
  studied (``Flipped'', which is defined in~Figure~\ref{fig:meshes}),
  $\dim \Spurious_h$ is strictly greater than $\sigma$.
\item[Stability] For all triangulations we have tested, the method
  seems at least reduced stable (i.e., stable after removal of
  spurious modes, if any), for $r = 2, 3$. This is in contrast to the
  situation for the Stokes equations, where for some triangulations,
  such as the diagonal triangulation, the method is not reduced stable for $r
  = 2, 3$, while for other triangulations, it is. For $r = 1$, reduced
  stability holds for some triangulations, but fails for others,
  including the diagonal triangulation.
\item[Convergence] We also studied convergence of the method on
  diagonal triangulations. For such meshes, the method was observed to
  be stable for $r > 1$, but unstable for $r = 1$. Theory predicts
  optimal convergence of $\pressure$ in $\Ltwo$ and $\velocity$ in
  $\Hdiv$ for a stable method and this is in fact what was observed.
  Such optimal convergence holds for $r = 2, 3, 4$, but not in the
  apparently unstable case $r = 1$. In the case $r \geq 4$, it is
  known that $\velocity$ converges at one order higher in $\Ltwo$ than
  in $\Hdiv$. No such increase of order was observed for $r < 4$.
\end{description}

The note is organized as follows. We introduce further notation and
summarize some key points of the theory of mixed finite element
methods in Section~\ref{sec:preliminaries}. Further, we derive some
eigenvalue problems associated with the stability conditions and give
a characterization of the Brezzi inf-sup constant for the mixed
Laplacian $\LaplaceInfsup$ in
Section~\ref{sec:eigenvalueproblems}. These eigenvalue problems
applied to the Stokes equations were also stated by
Malkus~\cite{art:Mal1981} and (in part) by Qin~\cite{thesis:Qin1994}
and provide a foundation for numerical investigations of the Brezzi
stability conditions. Based on these general results,
Section~\ref{sec:cPxdP} is devoted to the study of continuous
piecewise polynomials in two dimensions for the velocity and
discontinuous piecewise polynomials for the pressure.

\section{Notation and preliminaries}
\label{sec:preliminaries}

The notion of reduced stability of families of mixed finite element
spaces is a key point in this note. In order to make this notion
precise, this preliminary section aims to introduce notation and
summarize the stability notions for finite element discretizations of
abstract saddle point problems.

If $V$ is an inner product space, we denote the dual space by $V^{*}$,
the inner product on $V$ by $\langle \, \cdot \, , \, \cdot \,
\rangle_{V}$ and the induced norm by $\|\cdot\|_{V}$.  Let $\Omega$ be
an open and bounded domain in $\R^d$ with boundary $\partial
\Omega$. We let $H^m(\Omega)$, for $m = 0, 1, \dots$, denote the
standard Sobolev spaces of square integrable functions with $m$ weak
derivatives and denote their norm by $\|\cdot\|_m$. Accordingly,
$H^0(\Omega) = L^2(\Omega)$. The space of polynomials of degree $r$ on
$\Omega$ is denoted $\Poly{r}(\Omega)$. The space of vectors in $\R^d$
is denoted $\V$ and in general, $X(\Omega; \V)$ denotes the space of
vector fields on $\Omega$ for which each component is in
$X(\Omega)$. For brevity however, the space of vector fields in
$\LtwoVec$ with square integrable divergence is written $\Hdiv$ with
norm $\|\cdot\|_{\Div}$ and semi-norm $|\cdot|_{\Div} = \|\Div \cdot
\, \|_0$.  The subscripts and the reference to the domain $\Omega$
will be omitted when considered superfluous.

Let $\triang_h$ denote an admissible simplicial tessellation of
$\Omega$, $h$ measuring the mesh size of the tessellation. We shall
frequently refer to spaces of piecewise polynomials defined relative
to such, and label the spaces of continuous, and discontinuous,
piecewise polynomials of degree less than or equal to $r$ as follows.
\begin{align*}
  \cP{r} = \cP{r}(\triang_h) &= 
  \{ p \in H^1(\Omega)\, | \; p|_{K} \in \Poly{r}(K) 
  \quad \foralls K \in \triang_h \} \quad  r = 1, 2, \dots, \\
  \dP{r} = \dP{r}(\triang_h) &= 
  \{ p \in L^2(\Omega)\,|  \; p|_{K} \in \Poly{r}(K) 
  \quad \foralls K \in \triang_h \} \quad  r = 0, 1, \dots.
\end{align*}

The classical abstract saddle point problem reads as
follows~\cite{art:Bre1974, book:BreFor1991}: for given Hilbert spaces
$V$ and $Q$ and data $(f, g) \in V^{*} \times Q^{*}$, find
$(\velocity, \pressure) \in V \times Q$ satisfying
\begin{equation}
  \label{eq:saddlepoint}
  a(\velocity, \testvector) 
  + b(\testvector, \pressure)   + b(\velocity, \testscalar)
  =  \langle f, \testvector \rangle + \langle g, \testscalar \rangle
  \quad \forall \, (\testvector, \testscalar) \in V \times Q,
\end{equation}
where $a$ and $b$ are assumed to be continuous, bilinear forms on $V
\times V$ and $V \times Q$, respectively. We shall assume here and
throughout that $a$ is symmetric. Following~\cite{art:Bab1973}, there
exists a unique solution $(\velocity, \pressure)$
of~\eqref{eq:saddlepoint}, if and only if the continuous
Babu\v{s}ka inf-sup constant
\begin{equation}
  \label{eq:cont:Babuska}
  \contBabuska = 
  \inf_{0 \not = (\velocity, \pressure)}
  \sup_{0 \not = (\testvector, \testscalar)}
  \frac{a(\velocity, \testvector) + b(\testvector, \pressure) + b(\velocity, \testscalar)}
       {\|(\velocity, \pressure)\|_{\VV \times \QQ} 
         \|(\testvector, \testscalar)\|_{\VV \times \QQ}}
\end{equation}
is positive.  By~\cite{art:Bre1974}, this holds if and only if the continuous Brezzi
coercivity and Brezzi inf-sup constants are positive.  These are defined as
  \begin{align}
    \label{eq:cont:Coercivity:Infsup}
     \contCoercivity =
    \inf_{0 \not = \velocity \in \Kernel}
    \sup_{0 \not = \testvector \in \Kernel}
    \frac{a(\velocity, \testvector)}
         {\|\velocity\|_{\VV} \|\testvector\|_{\VV}} , \\
    \label{eq:cont:Infsup}
     \contInfsup =
    \inf_{0 \not = \testscalar \in \QQ}
    \sup_{0 \not = \testvector \in \VV}
    \frac{b(\testvector, \testscalar)}
         {\|\testvector\|_{\VV} \|\testscalar\|_{\QQ}},
  \end{align}
respectively,
where $\Kernel = \{\testvector \in \VV \, | \; b(\testvector, \testscalar)
= 0 \quad \foralls \testscalar \in \QQ\}$.

Given finite dimensional spaces $\VV_h \subset \VV$ and $\QQ_h \subset
Q$, defined relative to a tessellation $\triang_h$ of $\Omega$, the
Galerkin discretization of~\eqref{eq:saddlepoint} takes the form: Find
$(\velocity_h, \pressure_h) \in \VV_h \times \QQ_h$ satisfying
\begin{equation}
  \label{eq:disc:saddlepoint}
  a(\velocity_h, \testvector) 
  + b(\testvector, \pressure_h)
  + b(\velocity_h, \testscalar)  =
  \langle f, \testvector \rangle + \langle g, \testscalar \rangle
  \quad \forall \, (\testvector, \testscalar) \in \VV_h \times \QQ_h.
\end{equation}
On the discrete level, the 
Babu\v{s}ka inf-sup, Brezzi coercivity and Brezzi inf-sup
constants are defined as
\begin{align}
  \label{eq:Babuska}
  \Babuska &:=
  \inf_{0 \not = (\velocity, \pressure) \in \VV_h \times \QQ_h}
  \sup_{0 \not = (\testvector, \testscalar) \in \VV_h \times \QQ_h}
  \frac{a(\velocity, \testvector) + b(\testvector, \pressure) + b(\velocity, \testscalar)}
       {\|(\velocity, \pressure)\|_{\VV \times \QQ} \|(\testvector, \testscalar)\|_{\VV \times \QQ}},
\\
  \label{eq:Coercivity:Infsup}
  \Coercivity &:=
  \inf_{0 \not = \velocity \in \Kernel_h}
  \sup_{0 \not = \testvector \in \Kernel_h}
  \frac{a(\velocity, \testvector)}
       {\|\velocity\|_{\VV} \|\testvector\|_{\VV}}, \\
  \label{eq:Infsup}
   \Infsup &:=
  \inf_{0 \not =  \testscalar \in \QQ_h}
  \sup_{0 \not = \testvector \in \VV_h}
  \frac{b(\testvector, \testscalar)}
       {\|\testvector\|_{\VV} \|\testscalar\|_{\QQ}},
\end{align}
where
\begin{equation}
  \label{eq:Kernelh}
  \Kernel_h 
  = \{\testvector \in \VV_h\, |\; 
  b(\testvector, \testscalar) = 0 \quad \foralls \testscalar \in \QQ_h\}.
\end{equation}
For given $\VV_h \times \QQ_h$, there exists a unique
solution of~\eqref{eq:disc:saddlepoint} if and only if $\Coercivity$ and $\Infsup$
(or equivalently $\Babuska$) are positive. Furthermore, for a family of
discretization spaces $\VV_h \times \QQ_h$ parameterized by
$h$, if $\Coercivity$
and $\Infsup$ are uniformly bounded from below, then one obtains the quasi-optimal approximation
estimate~\cite{art:Bre1974}:
\begin{equation*}
  \|\velocity - \velocity_h\|_{\VV} 
  + \|\pressure - \pressure_h\|_{\QQ}
  \leq C \left (
  \inf_{\testvector \in \VV_h} \|\velocity - \testvector\|_{\VV} +
  \inf_{\testscalar \in \QQ_h} \|\pressure - \testscalar\|_{\QQ}
  \right ),
\end{equation*}
with $C$ depending only on the bounds for $\Coercivity$
and $\Infsup$ and bounds on the bilinear forms $a$ and $b$.
The uniform boundedness condition motivates the notion of stability
for pairs of finite element spaces.
\begin{definition}[Stable discretization]
A family of finite element discretizations $\{\VV_h \times \QQ_h\}_h$
is \emph{stable in $\VV \times \QQ$} if the Brezzi coercivity and
inf-sup constants $\Coercivity$ and $\Infsup$ (or equivalently the
Babu\v{s}ka inf-sup constant $\Babuska$) are bounded from below by a
positive constant independent of $h$.
\end{definition}
\noindent In accordance with standard terminology, we
say that $\{V_h \times Q_h\}_h$ satisfies the Brezzi coercivity
or inf-sup conditions if $\Coercivity$ or
$\Infsup$, respectively, are uniformly bounded from below.

There are families of discretizations that are not stable in the sense
defined above, but have a reduced stability property. More precisely,
for a pair $\VV_h \times \QQ_h$ consider the space of spurious modes
$\Spurious_h \subseteq \QQ_h$:
\begin{equation}
  \label{eq:kernelcomplement}
  \Spurious_h  
  = \{\testscalar \in \QQ_h \, | \; b(\testvector, \testscalar) = 0 \quad
  \foralls \testvector \in \VV_h\}.
\end{equation}
For a stable discretization, $\Spurious_h$ contains only the zero
element.  Indeed, $\Infsup = 0$ if and only if $\Spurious_h$ contains
non-zero elements. If $\Spurious_h$ is non-trivial, it is natural
to consider the reduced space $\Spurious_h^{\perp}$, the orthogonal
complement of $\Spurious_h$ in $Q_h$, in place of $\QQ_h$. This
motivates the definition of the reduced Brezzi inf-sup constant,
relating to the stability of $\VV_h \times \Spurious_h^{\perp}$:
\begin{equation}
  \label{eq:reduced:Infsup}
  \ReducedInfsup =
  \inf_{0 \not = \testscalar \in \Spurious_h^{\perp}}
  \sup_{0 \not = \testvector \in \VV_h}
  \frac{b(\testvector, \testscalar)}
       {\|\testvector\|_{\VV} \|\testscalar\|_{\QQ}},
\end{equation}
and the following definition of reduced stable. By definition,
$\ReducedInfsup \not = 0$.
\begin{definition}[Reduced stable discretization]
  A family of discretizations $\{\VV_h \times \QQ_h\}_h$ is
  \emph{reduced stable in $\VV \times \QQ$} if the Brezzi coercivity
  constant $\Coercivity$ and the reduced Brezzi inf-sup constant
  $\ReducedInfsup$, defined by \eqref{eq:Coercivity:Infsup} and
  \eqref{eq:reduced:Infsup}, are bounded from below by a positive
  constant independent of $h$.
\end{definition}

\section{Eigenvalue problems related to the Babu\v{s}ka-Brezzi constants}
\label{sec:eigenvalueproblems}

For a given set of discrete spaces, the Babu\v{s}ka and Brezzi
constants defined by~\eqref{eq:Babuska}--\eqref{eq:Infsup} can be
computed by means of eigenvalue problems. The form and properties of
the eigenvalue problem associated with the Brezzi inf-sup constant for
the Stokes equations were discussed by Qin
in~\cite{thesis:Qin1994}. Since also the Brezzi coercivity constant
plays a role for the mixed Laplacian, we begin this section by
deriving how the Brezzi coercivity constant can be computed by similar
eigenvalue problems. Actually, in our application in
Section~\ref{sec:cPxdP}, the Brezzi coercivity condition will be
automatic, however, we discuss it here in the abstract case, for the
sake of completeness. These eigenvalue problems were also stated, and
carefully analysed from an algebraic view-point, by
Malkus~\cite{art:Mal1981} in connection with the displacement-pressure
formulation of the linear elasticity equations. We continue by
observing that the continuous Brezzi inf-sup constant can be naturally
associated with the smallest eigenvalue of the Laplacian itself.

\subsection{Eigenvalue problems for the discrete Babu\v{s}ka-Brezzi constants}

Let $\VV_h \subset \VV$ and $\QQ_h \subset \QQ$ be given finite
dimensional spaces as before. It follows easily from the
definition that the Babu\v{s}ka inf-sup
constant $\Babuska = |\lambda_{\min}|$ when $\lambda_{\min}$ is the
smallest (in modulus) eigenvalue of the following generalized
eigenvalue problem: Find $\lambda \in \R$, $0 \not = (\velocity,
\pressure) \in \VV_h \times \QQ_h$ satisfying
\begin{equation}
  \label{eq:eigenvalue:Babuska}
  a(\velocity, \testvector) 
  + b(\testvector, \pressure) 
  + b(\velocity, \testscalar) 
  = \lambda 
  \left (\langle \velocity, \testvector \rangle_{\VV} 
  + \langle \pressure, \testscalar \rangle_{\QQ} \right )
  \quad \forall \, (\testvector, \testscalar) \in \VV_h \times \QQ_h.
\end{equation}
The following lemma identifies an eigenvalue problem
associated with the Brezzi inf-sup constant.
\begin{lem}[{Qin~\cite[Lemma 5.1.1 -- 5.1.2]{thesis:Qin1994}}]
  \label{lem:Infsup}
  Let $\lambda_{\min}$ be the smallest eigenvalue of the following
  generalized eigenvalue problem: Find $\lambda \in \R$, $0 \not =
  (\velocity, \pressure) \in \VV_h \times \QQ_h$ satisfying
  \begin{equation}
    \label{eq:eigenvalue:infsup}
    \langle \velocity, \testvector \rangle_{\VV}
    + b(\testvector, \pressure)  
    + b(\velocity, \testscalar) =
    - \lambda \langle \pressure, \testscalar \rangle_{\QQ}, 
    \quad \forall \, (\testvector, \testscalar) \in \VV_h \times \QQ_h.
  \end{equation}
  Then, $\lambda \geq 0$ and for $\Infsup$ defined
  by~\eqref{eq:Infsup}, $\Infsup = \sqrt{\lambda_{\min}}$.
\end{lem}
\noindent It can also be shown that the reduced Brezzi inf-sup
constant $\ReducedInfsup$ equals the square-root of the smallest
non-zero eigenvalue of~\eqref{eq:eigenvalue:infsup}~\cite[Theorem
5.1.1]{thesis:Qin1994}.

The Babu\v{s}ka and Brezzi inf-sup constants are thus easily computed,
given bases for the spaces $\VV_h$ and $\QQ_h$. As
for~\eqref{eq:eigenvalue:Babuska}, it is easily seen that the Brezzi
coercivity constant $\Coercivity = |\lambda_{\min}|$ where
$\lambda_{\min}$ is the smallest (in modulus) eigenvalue of the
eigenvalue problem: Find $\lambda \in \R$ and $0 \not = \velocity \in
\Kernel_h$ such that
\begin{equation}
  \label{eq:one}
  a(\velocity, \testvector) 
  = \lambda \langle \velocity, \testvector \rangle_{\VV}
  \quad \forall \; \testvector \in \Kernel_h.
\end{equation}
However, a basis for $\Kernel_h$ is usually not readily available,
thus hindering the actual computation of the eigenvalues
of~\eqref{eq:one}. Instead, the above eigenvalue problem over
$\Kernel_h$ can be extended to a generalized eigenvalue problem over
$\VV_h \times \QQ_h$: Find $\lambda \in \R$ and $ 0 \not = (\velocity,
\pressure) \in \VV_h \times \QQ_h$ such that
\begin{equation}
  \label{eq:two}
  a(\velocity, \testvector) 
  + b(\testvector, \pressure) + b(\velocity, \testscalar) 
  = \lambda \langle \velocity, \testvector \rangle_{\VV}
  \quad \foralls (\testvector, \testscalar) \in \VV_h \times \QQ_h.
\end{equation}
The following lemma establishes the equivalence between~\eqref{eq:one}
and~\eqref{eq:two}.
\begin{lem}
  \label{lem:equivalence_of_eigenvalues}
  If $(\lambda, \velocity)$ is an eigenpair of~\eqref{eq:one}, there
  exists a $\pressure \in \QQ_h$ such that $(\lambda, (\velocity,
  \pressure))$ is an eigenpair of~\eqref{eq:two}. Conversely, if
  $(\lambda, (\velocity, \pressure))$ is an eigenpair
  of~\eqref{eq:two} and $\velocity \not = 0$, then $\velocity \in
  \Kernel_h$ and $(\lambda, \velocity)$ is an eigenpair
  of~\eqref{eq:one}. For $0 \not = \pressure \in \Spurious_h$ and any
  scalar $\lambda$, $(\lambda, (0, \pressure))$ is an eigenpair
  of~\eqref{eq:two}, and these are the only eigenpairs
  of~\eqref{eq:two} with $\velocity = 0$.
\end{lem}
\begin{proof}
  Let $(\lambda, \velocity)$ be an eigenpair of~\eqref{eq:one}. Define
  $B_h: V_h \rightarrow \QQ_h$ such that $\langle B_h v, q
  \rangle_{\QQ} = b(v, q)$ for all $q \in \QQ_h$. Since $B_h :
  \Kernel_h^{\perp} \rightarrow B_h(\VV_h)$ is an isomorphism,
  $\pressure \in B_h(\VV_h) \subset \QQ_h$ is well-defined by
  \begin{equation*}
    \langle \pressure, \testscalar \rangle_{\QQ} = \lambda \langle
    \velocity, B_h^{-1} \testscalar \rangle_{\VV} - a(\velocity,
    B_h^{-1} \testscalar) \quad \foralls \testscalar \in B_h(\VV_h).
  \end{equation*}
  Then, for any $\testvector \in \Kernel_h^{\perp}$, $\pressure$
  satisfies
  \begin{equation*}
    b(\testvector, \pressure) = 
    \langle B_h \testvector, \pressure \rangle_{\QQ} = 
    \lambda \langle \velocity, \testvector \rangle_{\VV} 
    - a(\velocity, \testvector). 
  \end{equation*}
  Further, by definition $b(\testvector, \pressure) = 0$ for any
  $\testvector \in \Kernel_h$. Hence, by the assumption that
  $(\lambda, \velocity)$ is an eigenpair of~\eqref{eq:one}, $(\lambda,
  (\velocity, \pressure))$ satisfies~\eqref{eq:two}. The converse
  statement is obvious. Finally, letting $u = 0$ in~\eqref{eq:two}, we
  see that $(\lambda, (0, \pressure)$ satisfies~\eqref{eq:two} if and
  only if $\pressure \in \Spurious_h$, but for any $\lambda \in \R$.
\end{proof}
\noindent Note that, as a consequence of the last observation in
Lemma~\ref{lem:equivalence_of_eigenvalues}, if $\Spurious_h$ is
non-trivial, the generalized eigenvalue problem~\eqref{eq:two} is
computationally not well-posed since any scalar $\lambda$ is an
eigenvalue.

In the subsequent section, we shall numerically investigate the
stability of families of finite element discretizations $\VV_h
\times \QQ_h$ such that $\Div \VV_h \subseteq \QQ_h$ for the mixed
Laplacian, using the eigenvalue problem~\eqref{eq:eigenvalue:infsup}
in terms of standard bases for the spaces $\VV_h$ and $\QQ_h$. The
eigenvalue problem~\eqref{eq:two} does not enter, but would, were we
to investigate discretizations where $\Div \VV_h \not \subseteq
\QQ_h$.

\subsection{A characterization of the mixed Laplacian Brezzi inf-sup constant}
\label{subsec:characterization}

We now turn from the general setting to consider the $H(\Div) \times
L^2$ formulation of the mixed Laplacian \eqref{eq:laplace}.  In
Lemma~\ref{lem:characterization} below, we show that the Brezzi
inf-sup constant can be identified with the smallest eigenvalue of the
negative Laplacian. Consequently, if a discretization family $\{\VV_h
\times \QQ_h\}_h$ guarantees eigenvalue convergence for the mixed
Laplace eigenvalue problem and is such that $\Div \VV_h \subseteq
\QQ_h$, the Brezzi inf-sup constant of the discretization will
converge to the continuous Brezzi inf-sup constant.
\begin{lem}
  \label{lem:characterization}
  Let $\VV \subseteq \Hdiv$ and $\QQ \subseteq \Ltwo$ be such that
  that $\Div \VV \subseteq \QQ$. Consider the Brezzi inf-sup
  eigenvalue problem~\eqref{eq:eigenvalue:infsup} applied
  to~\eqref{eq:laplace}:
  \begin{equation}
    \label{eq:three}
    \langle \velocity, \testvector \rangle_{\Div}
    + \langle \Div \testvector, \pressure \rangle
    + \langle \Div \velocity, \testscalar \rangle 
    = - \lambda \langle \pressure, \testscalar \rangle 
    \quad \foralls
    (\testvector, \testscalar) \in \VV \times \QQ.
  \end{equation}
  Consider also the mixed Laplace eigenvalue problem:
  \begin{equation}
    \label{eq:four}
    \langle \hat \velocity, \testvector \rangle
    + \langle \Div \testvector, \hat \pressure \rangle
    + \langle \Div \hat \velocity, \testscalar \rangle 
    = - \hat \lambda
    \langle \hat \pressure, \testscalar \rangle 
    \quad \foralls
    (\testvector, \testscalar) \in \VV \times \QQ.
  \end{equation}
  Then, $(\lambda, (\velocity, \pressure))$ is an eigenpair
  of~\eqref{eq:three} if and only if $(\hat \lambda, (\hat \velocity,
  \hat \pressure))$ is an eigenpair of~\eqref{eq:four} where $\hat
  \lambda = \lambda (1 - \lambda)^{-1}$, $\hat \velocity = \velocity$
  and $\hat \pressure = (1 - \lambda) \pressure$. Moreover, in this
  case, $0 \leq \lambda < 1$ and $\hat \lambda \geq 0$. Also, $p \not
  = 0, \hat p \not = 0$.
\end{lem}
\begin{proof}
  Assume that $(\lambda, (\velocity, \pressure))$ is an eigenpair
  of~\eqref{eq:three}.  First, note that $\lambda \not = 1$. Letting
  $\lambda = 1$, $\testvector = \velocity$ and $q = - \Div \velocity$
  in~\eqref{eq:three}, implies that $\velocity = 0$. Further,
  $\testvector = 0$ and $\testscalar = \pressure$ gives that
  $\pressure = 0$. Hence, $\lambda = 1$ is only associated with the
  zero solution, which by definition, cannot form an eigenpair. Next,
  note that $\pressure \not = 0$, since otherwise implies that
  $\velocity = 0$, which again is impossible.  By the assumption $\Div
  \VV \subseteq \QQ$,
  \begin{equation}
    \label{eq:middle_calculation}
    \langle \Div \velocity, \Div \testvector \rangle 
    = - \lambda \langle \pressure, \Div \testvector \rangle
    \quad \foralls \testvector \in \VV.
  \end{equation}
  Taking $\testvector = \velocity$ in~\eqref{eq:middle_calculation}
  and letting $\testvector = \velocity$ and $\testscalar = (\lambda -
  1) \pressure$ in~\eqref{eq:three} show that $\|u\|^2 = \lambda(1 -
  \lambda) \|\pressure\|^2$. So, $0 \leq \lambda < 1$. The combination
  of~\eqref{eq:middle_calculation} and~\eqref{eq:three} gives
  \begin{equation*}
    \langle \velocity, \testvector \rangle
    +  \langle \Div \velocity, \testscalar \rangle 
    + (1 - \lambda) \langle \Div \testvector, \pressure \rangle
    = - \lambda \langle \pressure, \testscalar \rangle .
  \end{equation*}
  Finally, letting $\hat \velocity = \velocity$, $\hat{\pressure} = (1
  - \lambda) \pressure$ and $\hat \lambda = \lambda (1 -
  \lambda)^{-1}$, gives that $(\hat \lambda, (\hat \velocity, \hat
  \pressure))$ solves~\eqref{eq:four}. The converse holds by similar
  arguments.
\end{proof}

The equivalence demonstrated in the lemma above affords a simple
characterization of the Brezzi inf-sup constant for the mixed
Laplacian. The eigenvalue problem~\eqref{eq:three} with $\VV = \Hdiv$
and $\QQ = \Ltwo$ is the eigenvalue problem associated with the
continuous Brezzi inf-sup constant $\contLaplaceInfsup$,
cf.~Lemma~\ref{lem:Infsup}. Hence, $\contLaplaceInfsup$ is the
square-root of the smallest eigenvalue of~\eqref{eq:three}. On the
other hand, the eigenvalue problem~\eqref{eq:four} is a mixed weak
formulation of the standard eigenvalue problem for the negative
Laplacian with Dirichlet boundary conditions, given in strong form
below:
\begin{equation}
  \label{eq:laplace:eigenvalue}
  - \Delta \hat \pressure = 
  \hat \lambda \, \hat \pressure \text{ in } \Omega, \quad  
  \hat \pressure  = 0 \text{ on } \partial \Omega.
\end{equation}
Thus, if $\mu$ is the smallest eigenvalue
of~\eqref{eq:laplace:eigenvalue}, $\contInfsup^{\Div} =
\sqrt{\mu (1 + \mu)^{-1}}$.
\begin{remark}
An alternative eigenvalue problem arises from noting that, under the
assumption $\Div \VV \subseteq Q$,~\eqref{eq:four} with $\testvector =
0$ implies that $\hat \lambda \, \hat \pressure = - \Div \hat
\velocity$. Hence, if $(\hat \lambda, (\hat \velocity, \hat
\pressure))$ solves~\eqref{eq:four}, then either $\hat \lambda = 0,
\hat \velocity = 0$ and $\hat p \perp \Div \VV $, or $\hat \lambda >
0$, $\hat \velocity \not = 0$, and $(\hat \lambda, \hat \velocity)$ will be
an eigenpair of the problem:
\begin{equation*}
  \langle \Div \velocity, \Div \testvector \rangle = 
  \hat \lambda \langle \velocity, \testvector \rangle
  \quad \forall \, \testvector \in \VV.
\end{equation*}
This eigenproblem was studied in~\cite{art:BofDurGas1999}.
\end{remark}

Now, consider a stable discretization family $\VV_h \times
\QQ_h$ of~\eqref{eq:laplace} such that $\Div \VV_h
\subseteq \QQ_h$, with Brezzi inf-sup constants $\LaplaceInfsup$. Let $\mu_h$ denote the smallest eigenvalue approximation
of~\eqref{eq:four} by $\VV_h \times \QQ_h$. As a consequence of the
preceding considerations, if $\mu_h \rightarrow \mu$, then $\LaplaceInfsup \rightarrow
\contLaplaceInfsup$. In other words, if the discretization family is
stable, satisfies $\Div \VV_h \subseteq \QQ_h$, and gives eigenvalue
convergence, then the Brezzi inf-sup constant will converge to the
continuous Brezzi inf-sup constant. Note however, that the discrete
stability conditions are \emph{not} sufficient to ensure the convergence of
approximations to the eigenvalue
problem~\eqref{eq:four}~\cite{art:ArnFalWin2006, art:BofBreGas1999}.

Mixed finite element discretizations of~\eqref{eq:laplace} based on
the Raviart-Thomas~\cite{art:RavTho1977} and
Brezzi-Douglas-Marini~\cite{art:BreDouMar1985} families of $H(\Div)$
conforming elements are known to give eigenvalue convergence,
and hence $\LaplaceInfsup \rightarrow
\contInfsup^{\Div}$. Some cases where $\LaplaceInfsup$ seems to be
uniformly bounded in $h$, but $\LaplaceInfsup \not \rightarrow
\contInfsup^{\Div}$ are exemplified in the subsequent section. Finally,
note that if $\Omega$ is the unit square: $\Omega = (0, 1)^2$, the
smallest eigenvalue of~\eqref{eq:laplace:eigenvalue} is $2 \pi^2$ and
so
\begin{equation}
  \label{eq:exact:inf-sup}
  \contInfsup^{\Div} = \sqrt{\frac{2 \pi^2}{1 + 2 \pi^2}} 
  \approx 0.975593 .
\end{equation}

\section{Lower order Lagrange elements for the mixed Laplacian}
\label{sec:cPxdP}

From here on, we restrict our attention to finite element
discretizations of the mixed Laplacian~\eqref{eq:laplace} on a
polygonal domain $\Omega \subset \R^2$. The primary aim is to examine
the stability, or reduced stability, and convergence properties of
Lagrange elements, that is, continuous piecewise
polynomials for the vector variable and discontinuous piecewise
polynomials for the scalar variable:
\begin{equation}
  \label{eq:cPxdP}
  \VV_h \times \QQ_h = \cP{r}(\triang_h; \V) \times
  \dP{r-1}(\triang_h),
\end{equation}
for $r = 1, 2, \dots$.  Although the Brezzi conditions are in general
not satisfied for these discretizations, stability or reduced
stability may be identified on families of structured triangulations.
The pair~\eqref{eq:cPxdP} is clearly such that $\Div \VV_h \subseteq
\QQ_h$. Therefore, the stability of the discretization relies on a
uniform bound for the Brezzi inf-sup constant only. Further, a uniform
lower bound on the Brezzi inf-sup constant for the Stokes equations
induces the corresponding bound for the mixed Laplacian. Hence, the
results on the reduced stability of this element pair for the Stokes
equations can be directly applied to the mixed Laplacian. In the
following, new numerical evidence is presented and compared to the
known results.

The stability of the $\cP{r}(\V) \times \dP{r-1}$ family of elements,
for both the Stokes equations and the mixed Laplacian, depends on the
polynomial degree $r$ and the structure of the triangulation
$\triang_h$. For triangulations that have interior singular
vertices, the space of spurious modes $\Spurious_h$, defined
by~\eqref{eq:kernelcomplement} applied to~\eqref{eq:laplace}, will be
non-trivial. Here, an interior vertex is labeled singular if the
edges meeting at that vertex fall on two straight lines. Let $x$
denote an interior singular vertex and let $\omega_x$ be the star of
$x$. For any $r \geq 1$, there exists a $\pressure \in \Spurious_h$
such that $\pressure$ is supported in $\omega_x$~\cite{art:MorSco1975,
art:MorSco}. Consequently, letting $\sigma$ denote the number of
interior singular vertices of a triangulation, $\dim \Spurious_h \geq
\sigma$. Scott and Vogelius showed, for $r \geq 4$ that if there are
no interior singular vertices, then $\dim \Spurious_h = 0$ and so
$\LaplaceInfsup \ge \StokesInfsup > 0$. Moreover, they proved that for a
family of meshes without interior singular vertices, $\StokesInfsup$
remains bounded above zero as long as the meshes do not tend to
singularity as $h \rightarrow 0$. For the precise statement and more
details, see~\cite{art:ScoVog1985} or~\cite[Section
12.6]{book:BreSco2008}.

As we shall see below, for $r < 4$, the space of spurious modes may be
non-trivial even when there are no singular vertices. Further, for the
Stokes equations, more restrictive conditions than the above must be
placed on the triangulations in order to obtain a uniform bound for
the Stokes Brezzi inf-sup constant~\cite{thesis:Qin1994}. The
stability properties of these lower order discretizations for the
mixed Laplacian is the main question of interest in the following.

\begin{remark}
We shall not consider the pairing $\cP{r}(\V) \times \dP{s}$
except for $s=r-1$. This choice is easily motivated. First, a
dimension count shows that the pairing of continuous piecewise
polynomial vector fields with discontinuous piecewise polynomials of
the same or higher degree must have a non-trivial space of spurious
modes. Second, although the Brezzi inf-sup constant is uniformly
bounded for the pairs $\cP{r}(\V) \times \dP{r-2}$, $r = 2, 3, \dots$,
the Brezzi coercivity constant for the mixed Laplacian is not
uniformly bounded, and thus stability fails.
\end{remark}

\subsection{Stability}

\begin{figure}
  \begin{center}
    \subfigure[Diagonal]{
    \includegraphics[width=0.30\textwidth]{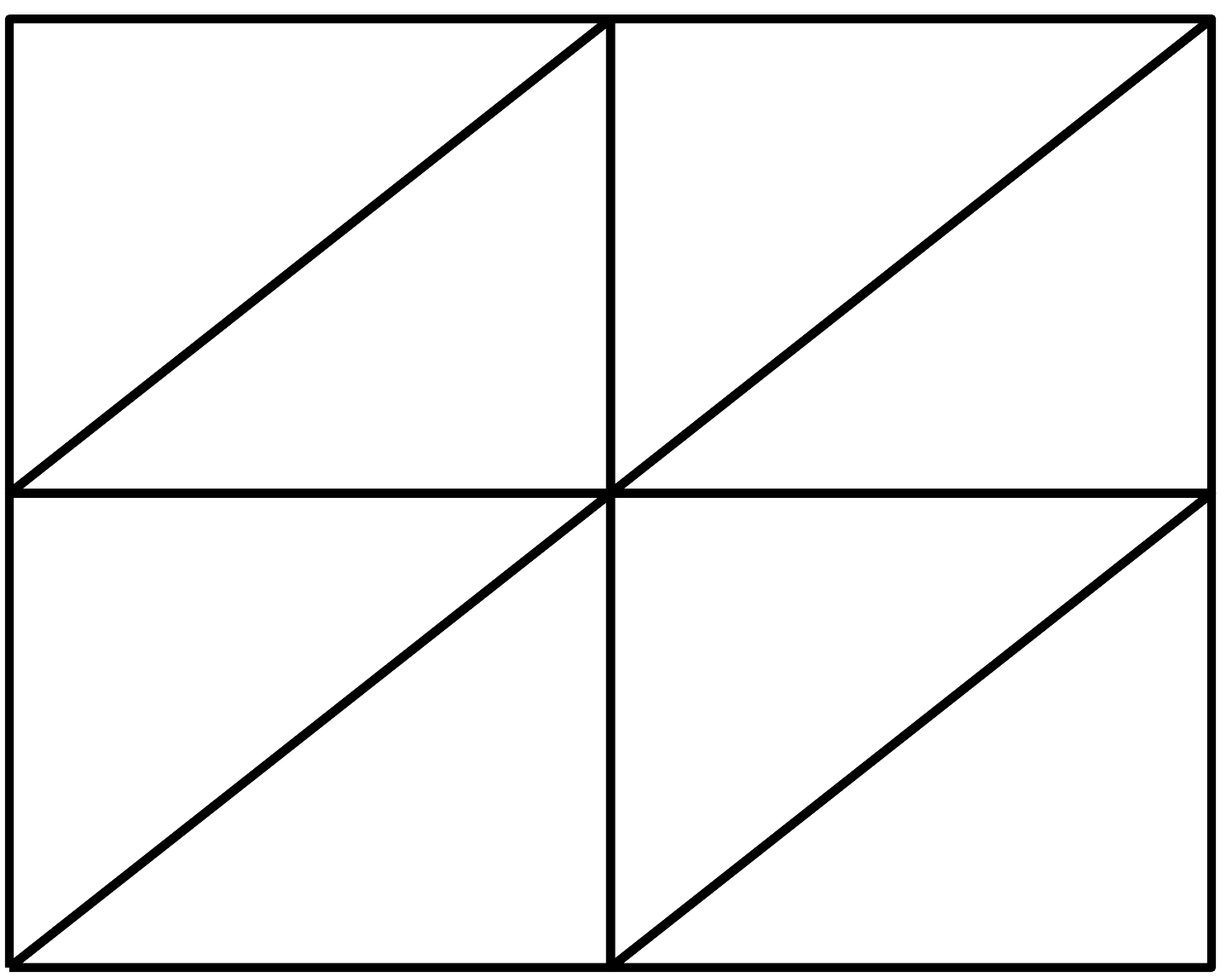}
    \label{mesh:diagonal}
    }
    \subfigure[Flipped]{
    \includegraphics[width=0.30\textwidth]{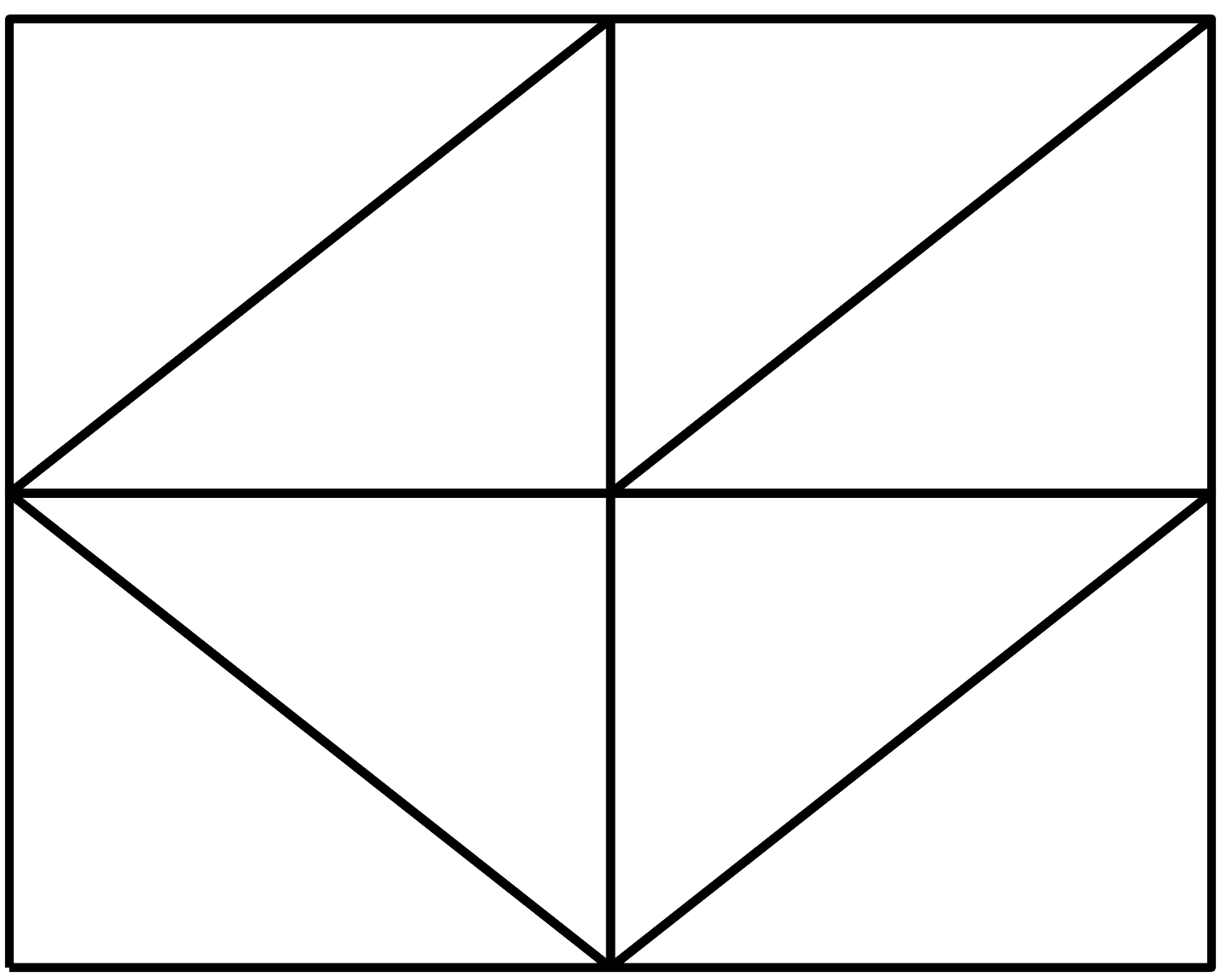} 
    \label{mesh:flipped}
    }
    \subfigure[Zigzag]{
    \includegraphics[width=0.30\textwidth]{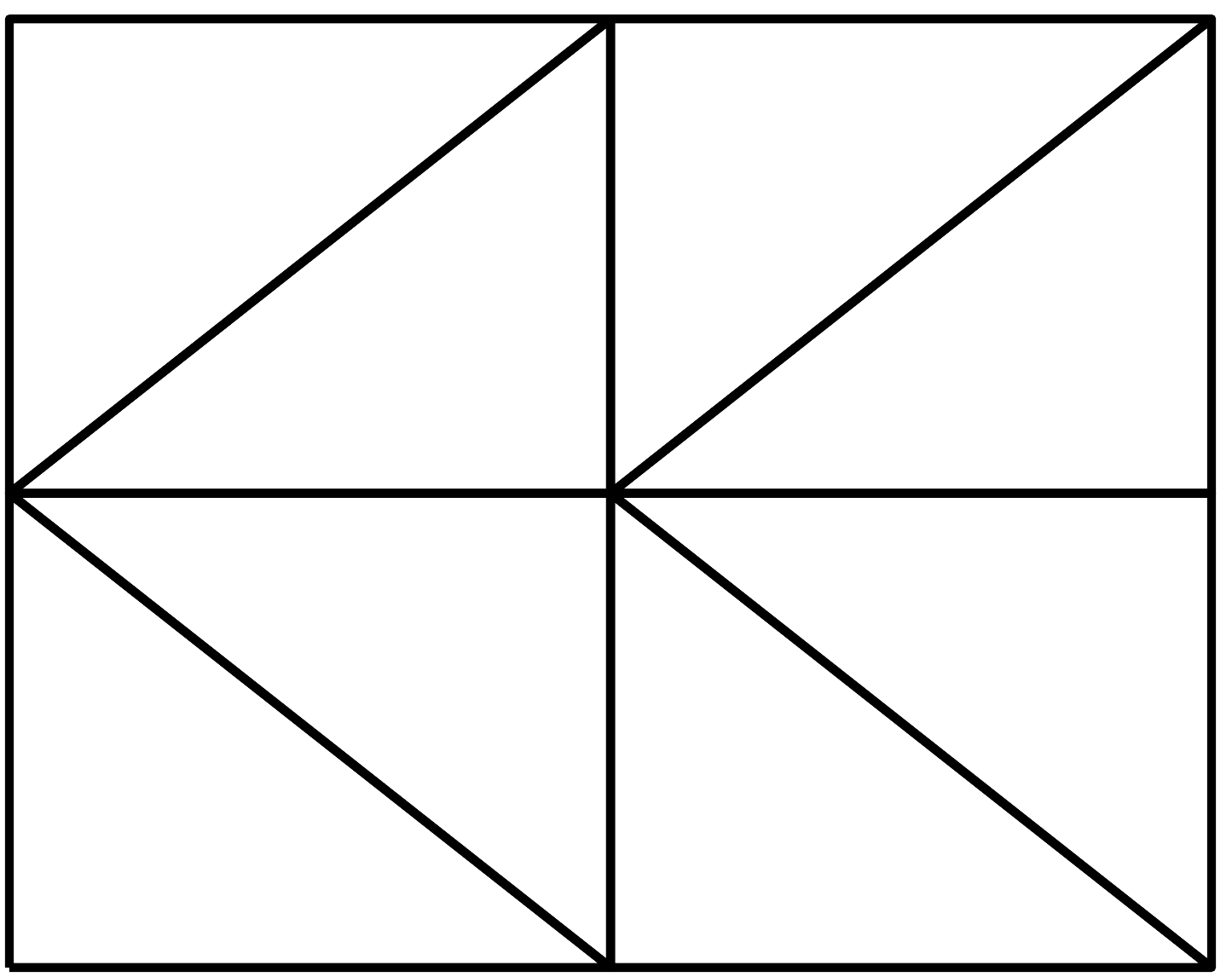} 
    \label{mesh:zigzag}
    }
    \subfigure[Crisscross]{
    \includegraphics[width=0.30\textwidth]{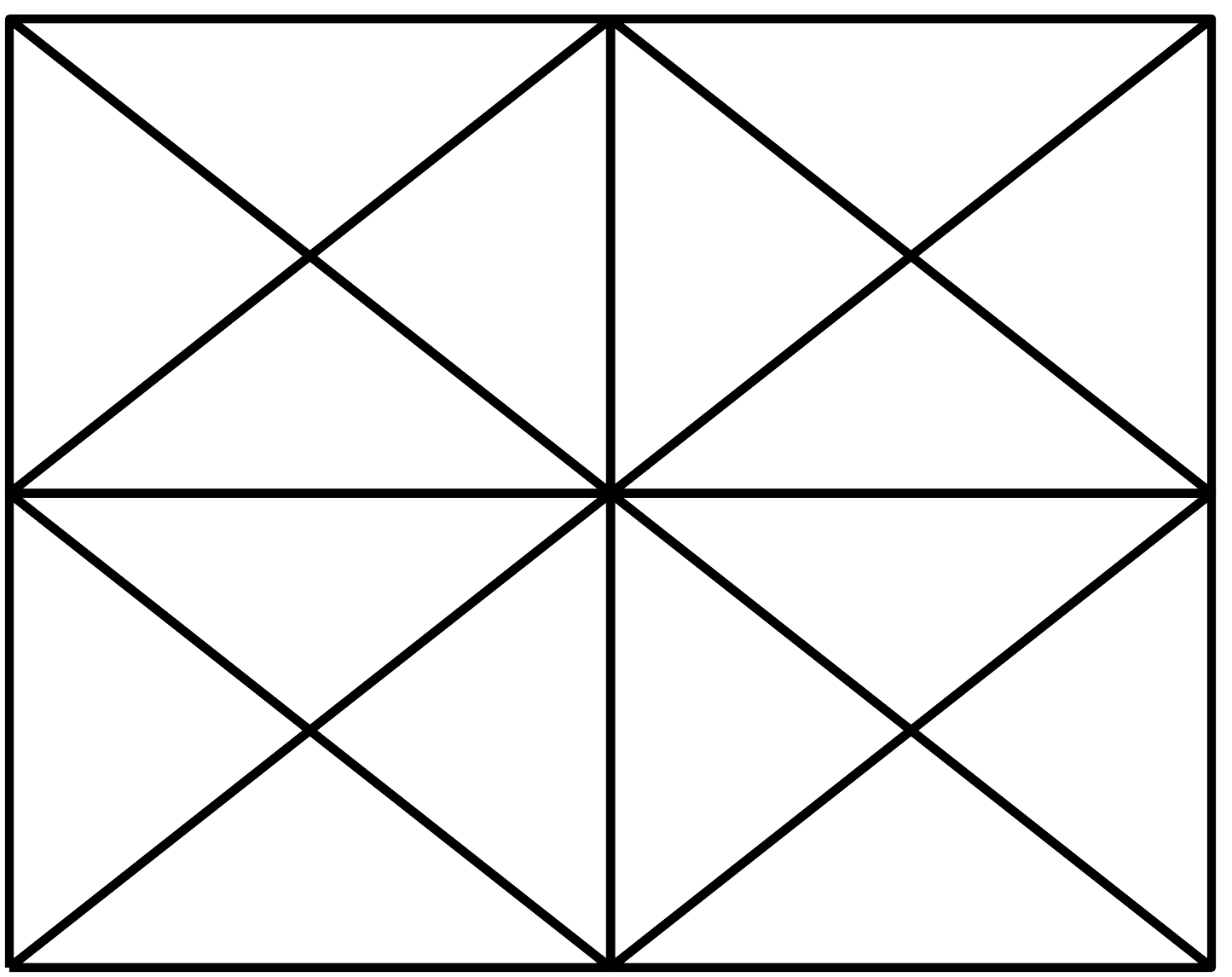}
    \label{mesh:crisscross}
    }
    \subfigure[Union Jack]{
    \includegraphics[width=0.30\textwidth]{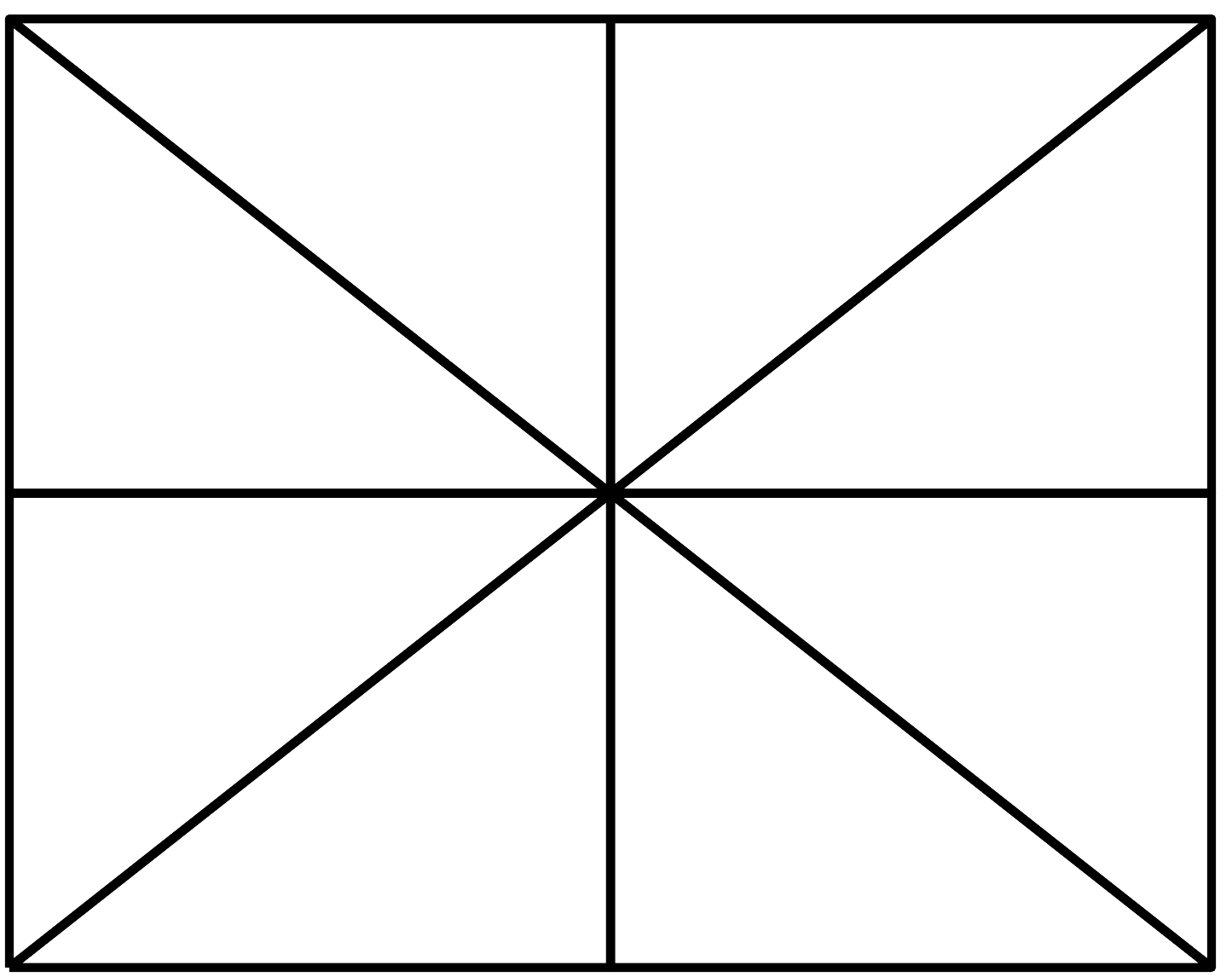}
    \label{mesh:unionjack}
    }
    \caption{Structured $2 \times 2$ triangulations of the unit
    square.}
    \label{fig:meshes} 
  \end{center}
\end{figure}

In the spirit of~\cite[Section 5]{thesis:Qin1994}, we aim to
numerically investigate the stability of $\cP{r}(\V) \times \dP{r-1}$
for $r = 1, 2, 3$ on certain families of structured triangulations of
the unit square.  The triangulation patterns considered are
illustrated and labeled in Figure~\ref{fig:meshes}. For $n$ even, an
$n \times n$ triangulation of each family is constructed by first
partitioning the domain into $n \times n$ squares, and subsequently
dividing each block of $2 \times 2$ squares into triangles by the
respective patterns. For instance, an $n \times n$ diagonal
triangulation is formed by dividing the unit square into $n \times n$
subsquares, and dividing each subsquare into triangles by the positive
diagonal. Throughout, we identify $h = 1/n$ and assume that $n > 2$.
Observe that the diagonal, flipped, and zigzag triangulations contain
no interior singular vertices, while the crisscross and the Union Jack
triangulation contain $n^2$ and $n(n-2)/2$ interior singular
vertices, respectively. This is summarized in the first row of
Table~\ref{tab:singular_vertices}.
\begin{table}
  \begin{center}
    \begin{tabular}{c|ccccc}
      \toprule
      & Diagonal & Flipped & Zigzag & Crisscross & Union Jack \\
      \midrule
      $\sigma$ & $0$ & $0$ & $0$ & $n^2$ & $n(n-2)/2$ \\
      $\dim \Spurious_h$ & $0$ &
      $\begin{cases}
        (n/2 - 1)^2 & r = 1 \\
        0 & r = 2, 3
        \end{cases}$
      & $0$ & $n^2$ & $n(n-2)/2$ \\
      \bottomrule
    \end{tabular}
  \end{center}
  \vspace{0.5cm}
  \caption{The number of interior singular vertices $\sigma$ and the
  dimension of the space of spurious modes, $\dim \Spurious_h$, for
  labeled families of $n \times n$ triangulations of the unit square,
  cf.~Figure~\ref{fig:meshes}. For the flipped, zigzag and Union Jack
  meshes, $\dim \Spurious_h$ is conjectural.}
  \label{tab:singular_vertices}
\end{table}

Recall that $\dim \Spurious_h \geq \sigma$ for $r \geq 1$ and equality
holds for $r \geq 4$. Qin proved that equality holds for $1 \leq r
\leq 3$ in the case of the diagonal and the crisscross meshes and
numerically observed equality for the flipped mesh for $r =
2$~\cite{thesis:Qin1994}. Our own experiments show that equality holds
for the zigzag and Union Jack meshes for $1 \leq r \leq 3$. Equality
also holds for the flipped mesh when $r = 2, 3$, but not for $r =
1$. These results are summarized in the second row of
Table~\ref{tab:singular_vertices}.

We continue by studying the behaviour of the Brezzi inf-sup constants
on the above triangulations. The cases $r = 2, 3$ are considered
first, but we will return to the case $r = 1$ below.  For the Stokes
equations, it is known that the diagonal and crisscross
triangulation families exhibit very different behaviour for $r = 2,
3$~\cite{thesis:Qin1994}. Namely, although there are non-trivial
spurious modes on the crisscross triangulation family, the reduced
Brezzi inf-sup constant is uniformly bounded. In contrast, for the
diagonal family, the Brezzi inf-sup constant decays as approximately
$\bigO(h)$.  As the discretization is reduced stable for the Stokes
equations on crisscross triangulations, it is also reduced stable for
the mixed Laplacian. A natural question becomes whether the lack of
stability on diagonal triangulations for the Stokes equations is also
present for the mixed Laplacian.

In view of Lemma~\ref{lem:Infsup}, we shall make an attempt at
answering this question through a set of numerical experiments. For a
given $r$ and a given $\triang_h$, the smallest, and smallest
non-zero, eigenvalue of~\eqref{eq:three} for $\VV = \cP{r}(\triang_h,
\V)$, $\QQ = \dP{r-1}(\triang_h)$ give the Brezzi inf-sup and reduced
Brezzi inf-sup constant. These eigenvalues for the triangulation
families considered, computed using LAPACK,
SLEPc~\cite{software:SLEPc} and DOLFIN~\cite{software:DOLFIN}, are
given for $r = 1, 2, 3$ in Tables~\ref{tab:cP1xdP0}, \ref{tab:cP2xdP1}
and~\ref{tab:cP3xdP2}.  For the purpose of identifying spurious modes,
eigenvalues below a threshold of $10^{-4}$ have been tabulated to
zero\footnote{Had a smaller threshold been chosen, some of the zero
eigenvalues associated to interior singular vertices would have been
missed for $r = 2$ on the Union Jack mesh of size $n = 6$.}.
\begin{table}
  \begin{center}
    \begin{tabular}{c|cc|cc}
      \toprule
      & \multicolumn{2}{c|}{$\LaplaceInfsup$} & 
      \multicolumn{2}{c}{$\ReducedInfsup^{\Div}$ ($\dim \Spurious_h$)} \\ 
      \midrule
      n & Diagonal & Zigzag 
      & Flipped
      & Union Jack \\
      \midrule
      4  & 0.847171  & 0.791967 & 0.945496 \digit(1)  & 0.976985 \digit\digit(4)\\
      6  & 0.716677  & 0.626865 & 0.945619 \digit(4)  & 0.976271 \digit(12) \\
      8  & 0.605576  & 0.505968 & 0.947850 \digit(9)  & 0.975985 \digit(24) \\
      10 & 0.517707  & 0.420180 & 0.946138 (16)  & 0.975847 \digit(40) \\
      12 & 0.449060  & 0.357720 & 0.944833 (25)  & 0.975770 \digit(60) \\
      14 & 0.394963  & 0.310731 & 0.943880 (36)  & 0.975724 \digit(84) \\
      16 & 0.351684  & 0.274303 & 0.943142 (49)  & 0.975693 (112) \\ 
      \bottomrule
    \end{tabular}
  \end{center}
  \vspace{0.5cm}
  \caption{The mixed Laplacian (reduced) Brezzi inf-sup constant for
    $\cP{1}(\triang_h; \V) \times \dP{0}(\triang_h)$ on labeled
    structured families of triangulations $\triang_h$. The dimension
    of the space of spurious modes in parenthesis if non-trivial.}
  \label{tab:cP1xdP0}
\end{table}
\begin{table}
  \begin{center}
    \begin{tabular}{c|ccc|c}
      \toprule
      & \multicolumn{3}{c|}{$\LaplaceInfsup$} & 
      $\ReducedInfsup^{\Div}$ ($\dim \Spurious_h$) \\ 
      \midrule
      n & Diagonal & Zigzag 
      & Flipped
      & Union Jack \\
      \midrule
      4 & 0.975627 & 0.955956  & 0.943790 & 0.975628 \digit(4) \\
      6 & 0.975600 & 0.952460  & 0.940480 & 0.975603 (12) \\
      8 & 0.975595 & 0.951384  & 0.938717 & 0.975595 (24) \\
      10 & 0.975594 & 0.950906 & 0.937684 & 0.975594 (40) \\
      12 & 0.975594 & 0.950638 & 0.936992 & 0.975593 (60) \\
      14 & 0.975593 & 0.950458 & & \\
      \bottomrule
    \end{tabular}
  \end{center}
  \vspace{0.5cm}
  \caption{The mixed Laplacian (reduced) Brezzi inf-sup constant for
    $\cP{2}(\triang_h; \V) \times \dP{1}(\triang_h)$ on labeled
    structured families of triangulations $\triang_h$. The dimension
    of the space of spurious modes in parenthesis if non-trivial.}
  \label{tab:cP2xdP1}
\end{table}
\begin{table}
  \begin{center}
    \begin{tabular}{c|ccc|c}
      \toprule
      & \multicolumn{3}{c|}{$\LaplaceInfsup$} & 
      $\ReducedInfsup^{\Div}$ ($\dim \Spurious_h$) \\ 
      \midrule
      n & Diagonal & Zigzag &  Flipped & Union Jack \\
      \midrule
      4 & 0.972244 & 0.975594 & 0.975594 & 0.975594 \digit(4) \\ 
      6 & 0.967304 & 0.975593 & 0.975593 & 0.975593 (12) \\
      8 & 0.964845 & 0.975593 & 0.975593 & 0.975593 (24) \\
      10 & 0.963412 &  &  &  \\
      12 & 0.962484 & & & \\
      \bottomrule
    \end{tabular}
  \end{center}
    \vspace{0.5cm}
    \caption{The mixed Laplacian (reduced) Brezzi inf-sup constant for
      $\cP{3}(\triang_h; \V) \times \dP{2}(\triang_h)$ on labeled
      structured families of triangulations $\triang_h$. The dimension
      of the space of spurious modes in parenthesis if non-trivial.}
  \label{tab:cP3xdP2}
\end{table}

For the diagonal meshes, the numerical experiments indicate that in
contrast to Stokes, the mixed Laplacian Brezzi inf-sup constants are
bounded from below for both $r = 2, 3$. For the flipped and zigzag
meshes, experiments give similar results. Neither exhibits any
spurious modes. Moreover, while the Brezzi inf-sup constant decays
approximately as $\bigO(h)$ for the Stokes
equations~\cite{thesis:Qin1994}, it appears to be uniformly bounded
for the mixed Laplacian. For the Union Jack family, the same
experiment gives $n(n-2)/2$ spurious modes, but the reduced
Brezzi inf-sup constant again seems to be uniformly bounded. In
summary for $r = 2, 3$, the $\cP{r}(\V) \times \dP{r-1}$ elements
appear to be at least reduced stable for all the families considered.

With the discussion in Section~\ref{subsec:characterization} in mind,
we also note that the Brezzi inf-sup constant converges to the exact
value, given by~\eqref{eq:exact:inf-sup}, for some, but not all, of
these meshes. For $r = 2$, the Brezzi inf-sup constant seems to
converge to the exact value on the diagonal meshes, but not on the
flipped or the zigzag meshes. The situation is the opposite for $r =
3$. There, the Brezzi inf-sup constant seems to converge to the exact
value on the zigzag and flipped meshes, but not for the diagonal
meshes.

The situation is different and more diverse in the lowest-order case:
$r = 1$. Boffi et al.~proved that $\cP{1}(\V) \times \dP{0}$ is in fact
reduced stable for the mixed Laplacian on crisscross
meshes~\cite{art:BofBreGas1999}. It is not reduced stable for
Stokes~\cite{thesis:Qin1994}.  However, the element pair does not seem
to be stable on diagonal meshes. The values in the first column of
Table~\ref{tab:cP1xdP0} indicate that the Brezzi inf-sup constant
decays approximately as $\bigO(h)$. The same is the case for the
zigzag meshes. For the Union Jack meshes, the situation is similar to
the crisscross case. That is, the number of singular modes match the
number of interior singular vertices and the reduced Brezzi inf-sup
constant appears to be bounded from below. Finally, the flipped meshes
display a surprising behaviour. There seem to be $(n/2 - 1)^2$
spurious modes, even though there are no singular vertices. This is
the only case where we have observed $\dim \Spurious_h >
\sigma$. However, the reduced Brezzi inf-sup constant appears to be
uniformly bounded.

\subsection{Convergence}

In the previous subsections, we have investigated the stability of the
$\cP{r}(\V) \times \dP{r-1}$ elements. Now, we proceed to examine the
convergence properties of these elements on the diagonal
meshes. Conjecturing that $\cP{r}(\V) \times \dP{r-1}$ is stable on
this mesh family for $r \geq 2$, in accordance with the numerical
evidence presented above, the standard theory gives the error estimate
\begin{equation}
  \label{eq:error:estimate}
  \| \velocity - \velocity_h \|_{\Div} 
  + \| \pressure - \pressure_h \|_{0}
   \leq C h^{r} \left ( 
   \|\velocity\|_{r+1} + \|\pressure\|_{r}
   \right ) .
\end{equation}
For $r \geq 4$, the $L^2$ error estimate for the velocity can be
improved~\cite[Theorem 12.4.9]{book:BreSco2008}, thus yielding:
\begin{equation}
  \label{eq:error:estimate:l2}
  \| \velocity - \velocity_h \|_{0} 
   \leq C h^{r+1} \|\velocity\|_{r+1}.
\end{equation}

\begin{figure}
  \begin{center}
    \subfigure[Normalized pressure errors in $\|\cdot\|_0$.]{
    \includegraphics[width=0.6\textwidth]{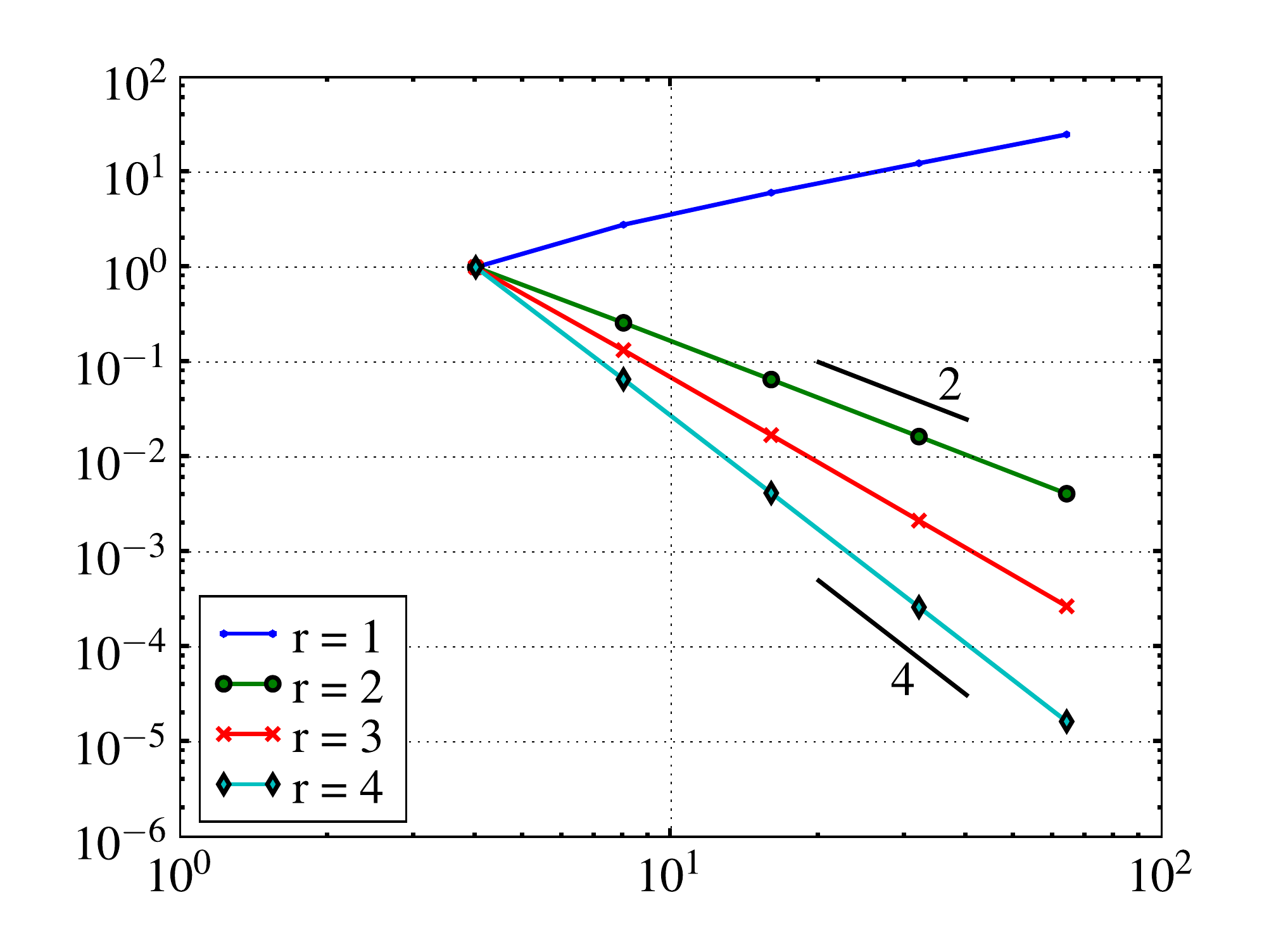}
    }
    \subfigure[Normalized velocity errors in $|\cdot|_{\Div}$]{
    \includegraphics[width=0.6\textwidth]{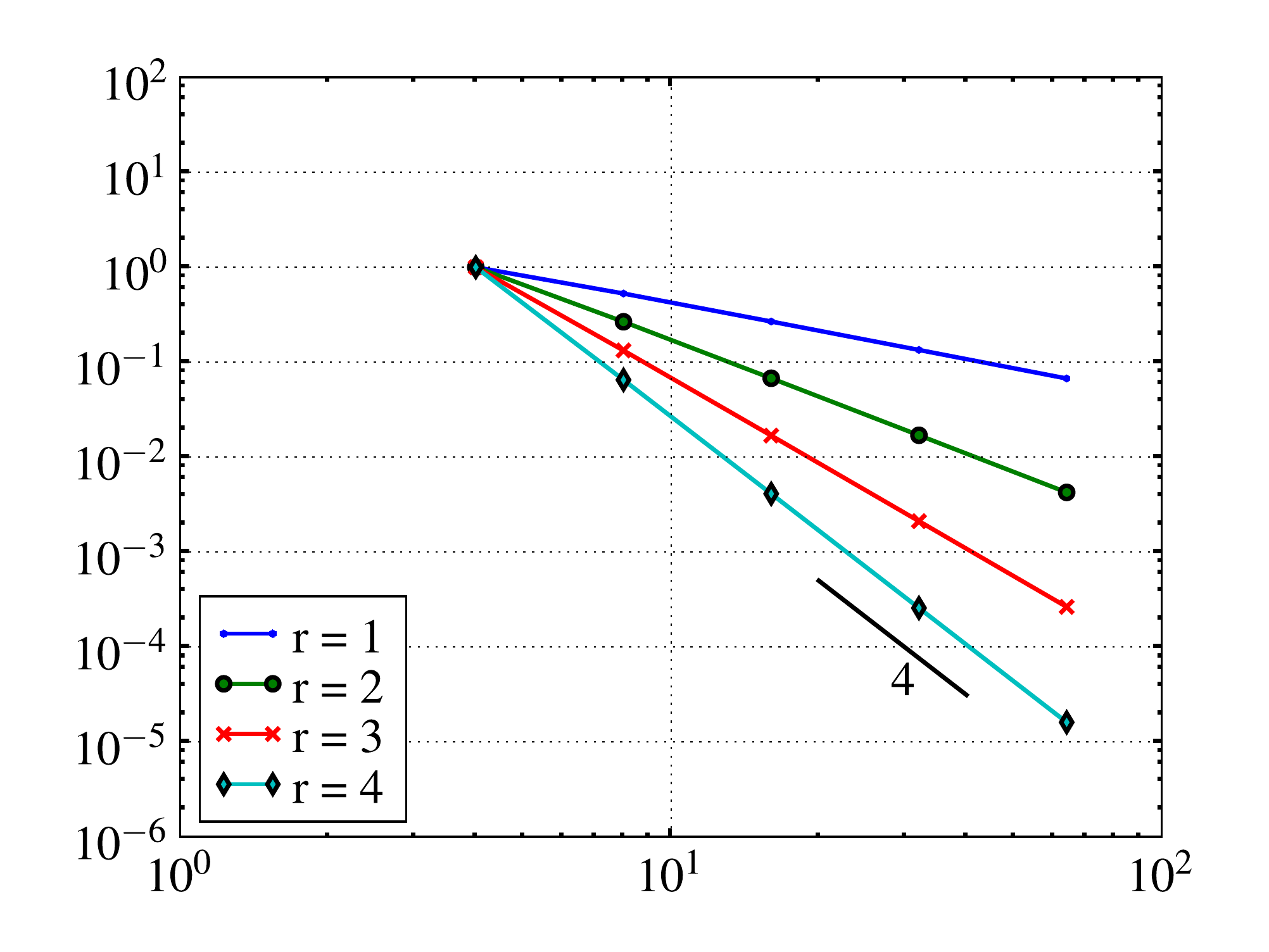}
    }
    \subfigure[Normalized velocity errors in $\|\cdot\|_0$]{
    \includegraphics[width=0.6\textwidth]{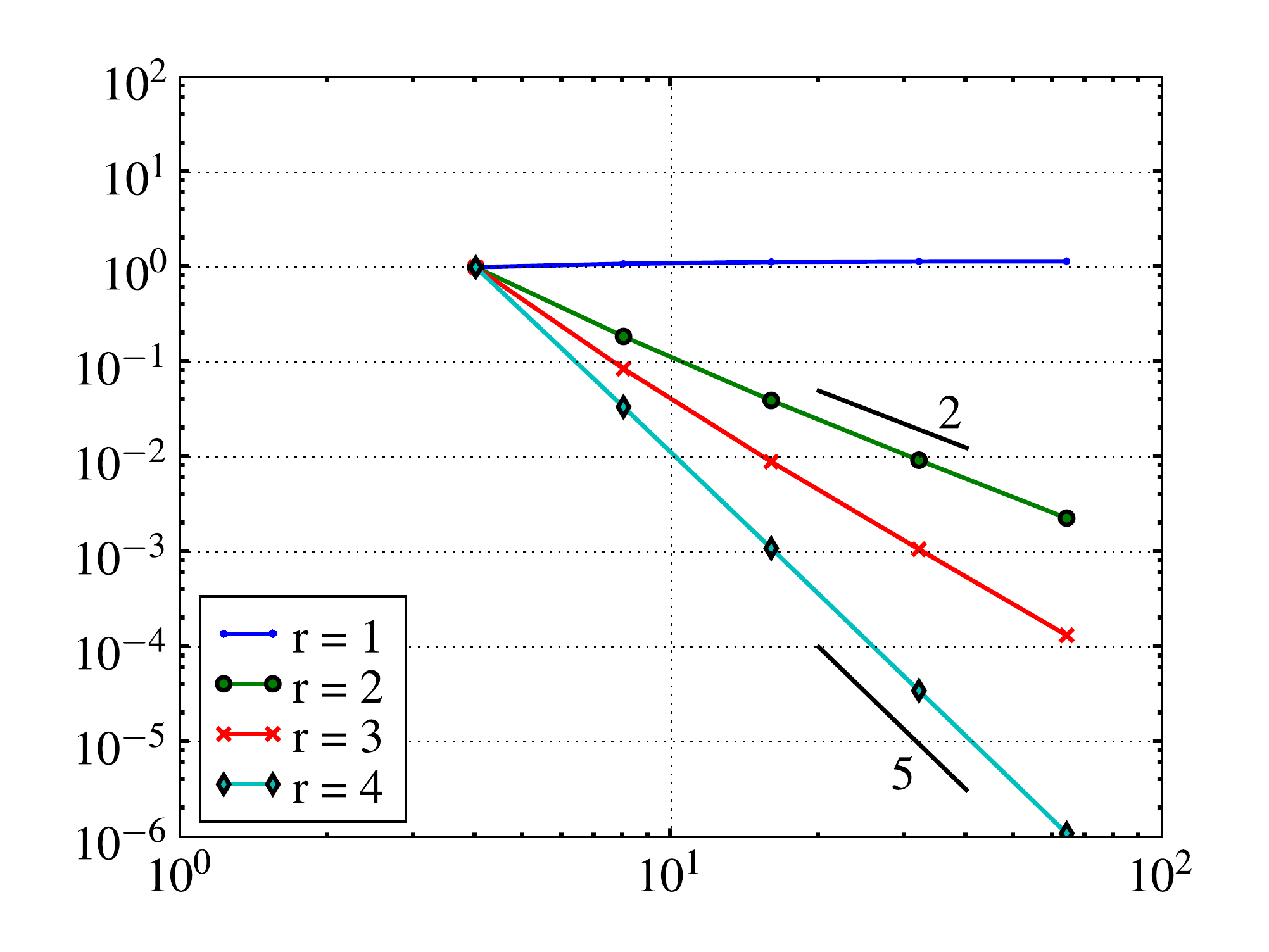}
    \label{fig:l2}
    }
    \caption{The errors of $\cP{r}(\V) \times \dP{r-1}$ approximations
    for $r = 1, 2, 3, 4$ on diagonal meshes versus mesh number
    $n$. The errors have been normalized, that is, multiplied by the
    inverse of the error at the smallest mesh $n = 4$. }
    \label{fig:convergence:cPxdP} 
  \end{center}
\end{figure}

In order to verify~\eqref{eq:error:estimate} and to see
whether~\eqref{eq:error:estimate:l2} appears to be attained for $r = 2, 3$,
we consider a standard smooth exact solution to the Laplacian with
pure Dirichlet boundary conditions:
\begin{equation}
  \label{eq:smooth_solution}
  \pressure(x,y) = \sin(2 \pi x) \sin(2 \pi y), \quad
  \velocity = \Grad \pressure, \quad
  g = \Div \velocity.
\end{equation}
The errors of the $\cP{r}(\V) \times \dP{r-1}$ approximations for $1
\leq r \leq 4$ on diagonal meshes can be examined in
Figure~\ref{fig:convergence:cPxdP}. To compute the errors, both the
source function $g$ and the exact solutions $u, p$ have been
represented by sixth order piecewise polynomial interpolants $\Pi_h g$
and $\Pi_h u, \Pi_h p$, whereupon the errors have been calculated
exactly (up to numerical precision).

For $r = 1$, we observed the discretization to be unstable on diagonal
meshes. As expected in this case, neither the pressure nor the
velocity approximation seems to converge in the $L^2$ norm. This
indicates that the estimates for the approximation error, based on the
standard estimates and the decaying Brezzi inf-sup constant, cannot be
improved. On the other hand, for $r > 1$, we observed the method to be
stable.  For $r = 2, 3, 4$, the orders of convergence in the $H(\Div)$
norm of the velocity and the $L^2$ norm of the pressure approximations
are indeed optimal, as predicted by~\eqref{eq:error:estimate}. 

The situation seems different for the convergence of the velocity
approximation in the $L^2$~norm. For $r \geq 4$, a convergence rate of
order $r+1$ is predicted by~\eqref{eq:error:estimate:l2}. This is also
observed for $r = 4$ in Figure~\ref{fig:l2}.  On the other hand, for
$r = 2$ and $3$ the rate of convergence appears to be of order $r$, and thus
one order suboptimal, indicating that the estimate~\eqref{eq:error:estimate:l2} does not
hold for $r<4$.

\bibliography{mxpoisson}
\bibliographystyle{abbrv} 

\end{document}